\documentclass[12pt]{amsart}
\usepackage{amsmath,amsthm,amssymb, amscd, amsfonts}

\theoremstyle{plain}
\newtheorem{theorem}{Theorem}[section]
\newtheorem{corollary}[theorem]{Corollary}

\newtheorem{proposition}[theorem]{Proposition}
\newtheorem{lemma}[theorem]{Lemma}
\newtheorem{question}[theorem]{Question}
\theoremstyle{definition}
\newtheorem{definition}[theorem]{Definition}

\newtheorem{example}[theorem]{Example}

\theoremstyle{remark}
\newtheorem{remark}[theorem]{Remark}

\numberwithin{equation}{section}\theoremstyle{plain}

\newcommand\id{\operatorname{id}}
\newcommand\ad{\operatorname{ad}}
\newcommand\Tr{\operatorname{Tr}}

\newcommand\co{\operatorname{co}}

\newcommand\cop{\operatorname{cop}}

\newcommand\op{\operatorname{op}}
\newcommand\Hom{\operatorname{Hom}}
\newcommand\Rep{\operatorname{Rep}}
\newcommand\rk{\operatorname{rk}}

\begin{document}
\title[$R$-matrices and Hopf algebra quotients]{$R$-matrices and Hopf algebra quotients}
\author{Sonia Natale}
\address{Facultad de Matem\'atica, Astronom\'\i a y F\'\i sica
\newline \indent
Universidad Nacional de C\'ordoba
\newline
\indent CIEM -- CONICET
\newline
\indent (5000) Ciudad Universitaria
\newline
\indent C\'ordoba, Argentina} \email{natale@mate.uncor.edu
\newline \indent \emph{URL:}\/ http://www.mate.uncor.edu/natale}

\thanks{This work was partially supported by CONICET, ANPCyT, SeCYT (UNC),
Fundaci\' on Antorchas and  FAMAF (Rep\'ublica Argentina).}


\date{\today.}

\begin{abstract} We study a natural construction of Hopf algebra quotients canonically associated to an
$R$-matrix in a finite dimensional Hopf algebra. We apply this
construction to show that a quasitriangular Hopf algebra whose
dimension is odd and square-free is necessarily cocommutative.
\end{abstract}

\maketitle

\section{Introduction and main results}

Quasitriangular Hopf algebras have been introduced by Drinfeld in
\cite{drinfeld}; these are Hopf algebras $H$ endowed with a
universal solution of the Quantum Yang-Baxter Equation $R \in H
\otimes H$, called an $R$-matrix. Such an $R$-matrix gives rise to
a map $\Phi_R: H^* \to H$, $$\Phi_R(f) : =
f(R^{(2)}r^{(1)})R^{(1)}r^{(2)}, \quad (R = r),$$  also studied by
Drinfeld in his paper \cite{Dracc}, and later interpreted in
categorical terms by S. Majid. Indeed, the \emph{transmutation}
map $\Phi_R$ turns out to be a morphism of \emph{braided} Hopf
algebras, after Majid's work.

Two extreme classes of quasitriangular Hopf algebras appear in
relation with the map $\Phi_R$:  \emph{triangular} Hopf algebras
are those for which $\Phi_R$ is trivial, while \emph{factorizable}
Hopf algebras are those for which $\Phi_R$ is an isomorphism.

Throughout this paper we shall work over an algebraically closed
field $k$ of characteristic zero.

Finite dimensional triangular Hopf algebras have been recently
classified by Etingof and Gelaki, based on results of Deligne on
symmetric categories \cite{eg-triangular}. On the other hand,
several important results are known to be true for factorizable
Hopf algebras \cite{re-se, eg, schneider}. In particular,
factorizable Hopf algebras are useful in the study of invariants
of $3$-manifolds. The Drinfeld double of a finite dimensional Hopf
algebra is always factorizable.

\smallbreak In this paper we give a construction of certain
canonical Hopf algebra quotients of a quasitriangular Hopf
algebra. Our construction is based on properties of transmutation
due to Majid \cite[Chapter 9.4]{majid}. Explicitly, the
construction goes as follows: any Hopf subalgebra of the dual is
'transmuted' into a normal coideal subalgebra, and this last one
is a sort of 'kernel' of a unique Hopf algebra quotient, in view
of the correspondence due to Takeuchi \cite{takeuchi}. This result
is contained in Theorem \ref{ccociente}. Several applications are
also given in Section \ref{estructura}.

An important special case of Theorem \ref{ccociente} is stated in
Theorem \ref{cociente}. It says that, in a sense, every
quasitriangular Hopf algebra is an 'extension' of the image of the
map $\Phi_R$ by a canonical triangular Hopf algebra.

\smallbreak We also discuss conditions, involving the $R$-matrix,
that guarantee normality of a quasitriangular Hopf algebra
quotient. In this context, we give a necessary and sufficient
condition for a map $H \to B$ to be normal; see Proposition
\ref{centraliz}. Combining these results with the classification
in \cite{aeg, eg2, eg-triangular} we prove in Theorem
\ref{hr-coconm} that if $(H, R)$ is a quasitriangular Hopf algebra
with $R \in kG(H) \otimes kG(H)$, then $H$ must be an extension of
a dual group algebra by a twisting of a modified supergroup
algebra.

\smallbreak We apply the results in Section \ref{estructura} to
prove a classification result for quasitriangular Hopf algebras.
The following open question appeared in \cite[Question
6.5]{andrusk}.

\begin{question}\label{pregunta} Let $N \geq 1$ be a square-free integer.
Are all Hopf algebras of dimension $N$ over $k$ necessarily
semisimple?
\end{question}

In the paper \cite{clspq} we proved that a quasitriangular Hopf
algebra of dimension $pq$, where $p$ and $q$ are odd prime numbers
is semisimple and therefore cocommutative. The following theorem
generalizes this result, and gives a partial answer to Question
\ref{pregunta} in the quasitriangular case.

\begin{theorem}\label{ap1} Let $(H, R)$ be a quasitriangular
Hopf algebra. Assume that $\dim H$ is odd and square-free. Then
$H$ is semisimple.

Moreover, $H \simeq kG$ where the group $G$ is a semidirect
product $G = \widehat \Gamma \rtimes F$, with $\Gamma$ cyclic. The
$R$-matrix corresponds under this isomorphism to an $F$-invariant
non-degenerate bicharacter $\langle \, , \, \rangle : \Gamma
\times \Gamma \to k^{\times}$ such that $\langle \, , \, \rangle
\langle \, , \, \rangle^{\op}$ is also non-degenerate.
\end{theorem}

Theorem \ref{ap1} will be proved in Section \ref{aplicaciones}.

\medbreak The paper is organized as follows. In Section \ref{fact}
we recall the definition and basic properties of quasitriangular
Hopf algebras. In Section \ref{apendice} we discuss the question
of normality of Hopf algebra quotients in the quasitriangular
case.  Section \ref{estructura} contains our main results on the
existence of Hopf algebra quotients attached to an $R$-matrix, and
some of its consequences. Finally, in Section \ref{aplicaciones}
we apply the results in previous sections to the classification of
quasitriangular Hopf algebras of square-free odd dimension.

\subsection*{Acknowledgement} The author thanks N. Andruskiewitsch
for discussions on the results of the paper \cite{aeg}.

\subsection{Preliminaries and notation}  Our
references for the theory of Hopf algebras are \cite{Mo, Sch}. For
a Hopf algebra $H$,  the antipode of $H$  will be denoted by
$\mathcal S$. The group of group-like elements in $H$ will be
denoted $G(H)$.

For an element $Z = \sum_i x_i \otimes y_i \in H \otimes H$, the
induced linear map $H^* \to H$, $p \to \sum_i \langle p, x_i
\rangle y_i$ will be denoted by $f_Z$. Thus, if $\dim H$ is
finite,  $Z \mapsto f_Z$ gives an isomorphism $H \otimes H \to
\Hom(H^*, H)$.

Also for such $Z$, the notation $Z_{21}$ will be used to indicate
the element $\tau(Z)$, where $\tau: H \otimes H \to H \otimes H$
is the usual flip map $\tau (x \otimes y) = y \otimes x$.

\medbreak Recall that Hopf algebra $H$ is called {\it simple} if
it contains no proper normal Hopf subalgebras in the sense of
\cite[3.4.1]{Mo}; $H$ is called {\it semisimple} (respectively,
{\it cosemisimple}) if it is semisimple as an algebra
(respectively, if it is cosemisimple as a coalgebra).

\section{Quasitriangular and factorizable Hopf algebras}\label{fact} Let $(H, R)$ be a
quasitriangular Hopf algebra over $k$. That is, $R \in H \otimes
H$ is an invertible element, called an \emph{$R$-matrix},
fulfilling the following conditions:

\begin{itemize}\item[(QT1)] $(\Delta \otimes \id)(R) = R_{13}R_{23}$.

\item[(QT2)] $(\epsilon \otimes \id)(R) = 1$.

\item[(QT3)] $(\id \otimes \Delta)(R) = R_{13}R_{12}$.

\item[(QT4)] $(\id \otimes \epsilon)(R) = 1$.

\item[(QT5)] $\Delta^{\cop}(h) = R \Delta(h) R^{-1}$, $\forall h \in H$.
\end{itemize}

The following relations with the antipode of $H$ are well-known:
\begin{equation}\label{inv-r} (\mathcal S \otimes \id)(R) = R^{-1}
= (\id \otimes \mathcal S^{-1})(R), \quad (\mathcal S \otimes
\mathcal S)(R) = R.\end{equation}

\medbreak There are Hopf algebra maps $f_R: {H^*}^{\cop} \to H$
and $f_{R_{21}}:H^* \to H^{\op}$ given, respectively, by
$$f_R(p) = \langle p, R^{(1)} \rangle R^{(2)},
\quad f_{R_{21}}(p) = \langle p, R^{(2)} \rangle R^{(1)},$$ for
all $p \in H^*$. Here, and elsewhere, we use the shorthand
notation $R = R^{(1)} \otimes R^{(2)} \in H \otimes H$. The images
of these maps will be indicated by $H_+$ and $H_-$, respectively.

Thus $H_+, H_-$ are finite dimensional Hopf subalgebras of $H$ and
$H_+ \simeq (H_-^*)^{\cop}$. The \emph{rank} of $R$ is defined as
$\rk R : = \dim H_+ = \dim H_-$.

We shall denote by $H_R = H_-H_+ = H_+H_-$ the minimal
quasitriangular Hopf subalgebra of $H$. See \cite{R}.

\begin{example}\label{cocomm} (See \cite{davydov}). Let $H = kG$ be a cocommutative Hopf algebra,
where $G = G(H)$ is a finite group. Suppose $(H, R)$ is a
quasitriangular structure on $H$. Since ${H^*}^{\cop}$ and $H^*$
are commutative, then $H_+$ and $H_- \simeq {H_+}^{*\cop}$ are
commutative and cocommutative Hopf subalgebras of $H$. In
particular, the Drinfeld double $D(H_+)$ is a commutative Hopf
algebra.

Write $H_R = k\Gamma$, where $\Gamma \subseteq G$ is a subgroup.
In view of \cite{R}, there is a Hopf algebra projection $D(H_+)
\to H_R = k\Gamma$, implying that $\Gamma$ is an \emph{abelian}
subgroup.

The set $\{e_a \}_{a \in \widehat \Gamma}$ is a complete set of
primitive orthogonal idempotents in $k\Gamma$, where $e_a =
\dfrac{1}{|\Gamma|}\sum_{g \in \Gamma}a(g^{-1})g$, $a \in \widehat
\Gamma$. Then $R$-matrix $R \in k\Gamma \otimes k\Gamma$ can be
written as
\begin{equation}\label{r-cocom}R = \sum_{a, b \in \widehat \Gamma} \rho(a, b) e_a \otimes
e_b,\end{equation} for a unique bilinear form $\rho: \widehat
\Gamma \times \widehat \Gamma \to k^{\times}$. Now Condition (QT5)
implies that $\Gamma$ is a \emph{normal} subgroup of $G$, and
moreover $\rho: \widehat \Gamma \times \widehat \Gamma \to
k^{\times}$ is an $\ad G$-invariant bilinear form on $\widehat
\Gamma$.

Conversely, it is not difficult to see that every pair $(\Gamma,
\rho)$, with $\Gamma \subseteq G$ a normal abelian subgroup and
$\rho: \widehat \Gamma \times \widehat \Gamma \to k^{\times}$ an
$\ad G$-invariant bilinear form, determines a unique
quasitriangular structure on $H = kG$ via the formula
\eqref{r-cocom}.

\end{example}

\medbreak Recall from \cite{re-se} that the quasitriangular Hopf
algebra $(H, R)$ is called {\it factorizable} if the map $\Phi_R :
H^* \to H$ is an isomorphism, where
\begin{equation}\label{phir} \Phi_R (p) = \langle p, Q^{(1)} \rangle Q^{(2)}, \quad p \in H^*;
\end{equation}
here, $Q = Q^{(1)} \otimes Q^{(2)} = R_{21} R \in H \otimes H$. In
other words, $\Phi_R = f_Q$.

\begin{remark}\label{rightspan} The image $\Phi_R(H^*)$ coincides with the subspace
of $H$ spanned by the right tensorands of $Q \in H \otimes H$.

Note that $\Phi_R (H^*) \subseteq H_R$. So every factorizable Hopf
algebra is minimal quasitriangular in the sense of \cite{R}.
\end{remark}

If on the other hand $\Phi_R = \epsilon 1$ (or equivalently,
$R_{21}R = 1 \otimes 1$), then $(H, R)$ is called
\emph{triangular}. Finite dimensional triangular Hopf algebras
were completely classified in \cite{eg-triangular}.

\medbreak We may also consider the map ${}_R\Phi = f_{Q_{21}} :
H^* \to H$, given by
\begin{equation} {}_R\Phi (p) = \langle p, Q^{(2)} \rangle Q^{(1)}, \quad p \in H^*.
\end{equation}

\begin{lemma}\label{left-right} ${}_R\Phi = \mathcal S  \Phi_R \mathcal S$.
In particular, ${}_R\Phi(H^*) = \mathcal S  \Phi_R(H^*)$.
\end{lemma}

\begin{proof}  We have $(\mathcal S \otimes \mathcal S) (R) = R$, and therefore
also $(\mathcal S \otimes \mathcal S) (R_{21}) = R_{21}$. Denote
as before $Q = R_{21} R$, then
$$(\mathcal S \otimes \mathcal S) (Q) =
(\mathcal S \otimes \mathcal S) (R) (\mathcal S \otimes \mathcal
S) (R_{21}) = R R_{21} = Q_{21}.$$ Therefore, for all $p \in H^*$,
$$ {}_R\Phi (p)
= \langle p, \mathcal S (Q^{(1)}) \rangle \mathcal S (Q^{(2)})
 = \mathcal S ( \langle \mathcal S (p), Q^{(1)} \rangle Q^{(2)}) =
\mathcal S (\Phi_R (\mathcal S (p) ) ).$$ \end{proof}

\begin{remark} Suppose that $H$ is semisimple. Let $R(H) \subseteq H^*$ and $Z(H) \subseteq H$ denote,
respectively, the character algebra and the center of $H$; so that
$R(H)$ coincides with the subalgebra of cocommutative elements in
$H^*$. Then $\Phi_R : R(H) \to Z(H)$ is an algebra map
\cite{Dracc}.

Note that for all $p, q \in R(H)$, $\langle p \otimes q, R_{21}R
\rangle = \langle q \otimes p, R_{21}R \rangle$; and therefore,
$\langle p, \Phi_R(q) \rangle = \langle q, \Phi_R(p) \rangle$.

So that the bilinear form $[\ , \ ]: R(H) \otimes R(H) \to k$,
$[p, q] = \langle p \otimes q, R_{21}R \rangle$ is symmetric.
Moreover, $(H, R)$ is factorizable if and only if $[\ , \ ]$ is
non-degenerate.
\end{remark}

\medbreak Let $\widehat H = \{ \chi_0 = \epsilon, \dots, \chi_n
\}$ be the set of irreducible characters of $H$. Let
$$s_{ij} : = \langle \chi_i, \Phi_R(\chi_j) \rangle = \langle \chi_j,
\Phi_R(\chi_i) \rangle, \quad 0 \leq i, j \leq n.$$ Then $(H, R)$
is factorizable if and only if the matrix $S = (s_{ij})_{0 \leq i,
j \leq n}$ is non-degenerate, if and only if the category $H$-mod
of finite dimensional $H$-modules is {\it modular}.

\begin{example}\label{doble} Recall that for a finite dimensional Hopf algebra $A$, its Drinfeld
double is a quasitriangular Hopf algebra. $D(A) = A^{* \cop}
\otimes A$ as a coalgebra, with a canonical $R$-matrix $\mathcal R
= \sum_i a^i \otimes a_i$, where $(a_i)_i$ is a basis of $A$ and
$(a^i)_i$ is the dual basis.

It is well-known that the Drinfeld double $(D(A), \mathcal R)$ is
factorizable \cite{re-se}. We have $D(A)_+ = A$, $D(A)_- =
{A^*}^{\cop}$.
\end{example}

\medbreak Our next result relates factorizability with
semisimplicity of a quasitriangular Hopf algebra. We showed in
\cite[Corollary 2.5]{clspq} that if $(H, R)$  is a factorizable
odd-dimensional Hopf algebra such that all proper Hopf subalgebras
of $H$ are semisimple, then $H$ is itself semisimple.

The following lemma gives a further criterion in this direction.
It will be used in Section \ref{aplicaciones}.

\begin{lemma}\label{fact-odd} Let $(H, R)$ be a factorizable Hopf algebra. Assume
in addition that $\dim H$ is odd and $\rk R = \dim H$. Then $H$ is
semisimple. \end{lemma}

\begin{proof} Recall that $H$ is semisimple if and only if $\Tr \mathcal S^2 \neq
0$ \cite{radford-tr}. Therefore if  $\mathcal S^4 = \id$ and $\dim
H$ is odd, then $H$ must be semisimple.

Let $g \in G(H)$, $\alpha \in G(H^*)$ be the modular elements of
$H$. By \cite{radford}, we have $\alpha = 1$. On the other hand,
the assumption on the rank of $R$ says that $f_R$ is an
isomorphism. Therefore $f_{R_{21}}(\alpha) = g^{-1}$ by
\cite[Corollary 2.10]{gelaki}. This implies that also $g = 1$.
Then $\mathcal S^4 = \id$ \cite[Proposition 10.1.14]{Mo}. Since
$\dim H$ is odd, this implies that $H$ is semisimple as claimed.
\end{proof}

\section{Normal quotients of quasitriangular Hopf
algebras}\label{apendice}

Let $H$ and $B$ be  finite dimensional Hopf algebras over $k$ and
let $\pi: H \to B$ be a surjective Hopf algebra map.

Then $H^{\co \pi}: = \{ h \in H: (\id \otimes \pi) \Delta (h) = h
\otimes 1\}$, is a left coideal subalgebra of $H$. Similarly,
${}^{\co \pi}H: = \{ h \in H: (\pi \otimes \id) \Delta (h) = 1
\otimes h\}$ is a right coideal subalgebra of $H$. We shall use
the notation $H^{\co B} := H^{\co \pi}$, ${}^{\co B}H := {}^{\co
\pi}H$, when no ambiguity arises.

Moreover, $H^{\co \pi}$ is stable under the left adjoint action of
$H$, defined as
$$\ad_h(a): = h_1 a \mathcal S(h_2),$$ for all $h, a \in H$. A
subspace of $H$ stable under the adjoint action will be called
\emph{normal}.

\begin{remark}\label{lcs-rcs} Consider the left (respectively, right) action of
$B^*$ on $H$ given by $f \rightharpoonup h = \langle f, \pi h_2
\rangle h_1$, (respectively, $h \leftharpoonup f = \langle f, \pi
h_1 \rangle h_2$). Then we have $H^{\co \pi} = {}^{B^*}H$ and
${}^{\co \pi}H = H^{B^*}$.
\end{remark}

\begin{definition} We shall say that a sequence $K \subseteq H \overset{\pi}\to
B$ is \emph{exact} if $\pi$ is surjective Hopf algebra map and $K
= H^{\co B}$ or $K = {}^{\co B}H$.
\end{definition}

\medbreak The surjection $\pi$ is called {\it normal} if $H^{\co
\pi} = {}^{\co \pi}H$. In this case $K = H^{\co \pi}$ is a
(normal) Hopf subalgebra of $H$ and there is an exact sequence of
Hopf algebras $1 \to K \to H \to B \to 1$. Furthermore, $B \simeq
H/HK^+$ as Hopf algebras.

\begin{lemma} The following are
equivalent:

\begin{itemize}\item[(i)] The surjection $\pi$ is normal.
\item[(ii)] $H^{\co \pi}$ is a right coideal of $H$.
\item[(iii)] $H^{\co \pi}$ is a subcoalgebra of $H$.
\end{itemize}
\end{lemma}

\begin{proof} We omit the details of the proof.
\end{proof}

\begin{remark}\label{sub-hopf} Suppose that $A$ is a Hopf subalgebra of $H$.
Then  $A \subseteq H^{\co \pi}$ if and only if $A \subseteq
{}^{\co \pi}H$ if and only if $\pi\vert_A = \epsilon$.
\end{remark}

\bigbreak \emph{From now on $(H, R)$ will be a finite dimensional
quasitriangular Hopf algebra and $q: H \to B$ will be a surjective
Hopf algebra map.}

\medbreak We aim to give a necessary and sufficient condition for
the map $q$ to be normal. We point out that a result in this
direction has been obtained by Masuoka in \cite[Corollary
5]{masuoka}, where a sufficient condition for normality of a
special kind of quotient is given. Contrary to our approach, this
condition is independent of the $R$-matrix involved.

\medbreak Recall that $H_+ \subseteq H$ denotes the Hopf
subalgebra of $H$ spanned by the right tensorands of $R$.

\begin{proposition}\label{cond-norm} Suppose $R \in H \otimes H^{\co B}$. Then $q: H \to B$ is
normal. \end{proposition}

\begin{proof} We need to show that $H^{\co B}$ is
a subcoalgebra, hence a Hopf subalgebra, of $H$.

We know that $\Delta(H^{\co B}) \subseteq H \otimes H^{\co B}$.
The assumption $R \in H \otimes H^{\co B}$ implies that also
$R^{-1} = (\mathcal S \otimes \id)(R) \in H \otimes H^{\co B}$.
Let $h \in H^{\co B}$, so that $\Delta(h) \in H \otimes H^{\co
B}$. On the other hand, we have
$$\Delta^{\cop}(h) = R \Delta(h)R^{-1} \in H \otimes H^{\co B}.$$
This implies that $\Delta(h) \in (H \otimes H^{\co B})  \cap
(H^{\co B} \otimes H) = H^{\co B} \otimes H^{\co B}$. Therefore
$H^{\co B}$ is a subcoalgebra of $H$ as claimed. \end{proof}

\begin{corollary}\label{inclusion} Suppose that $H_{+} \subseteq H^{\co B}$. Then
$q: H \to B$ is normal. \end{corollary}

\begin{proof} $H_{+}$ is the smallest Hopf subalgebra of $H$ with $R \in H \otimes H_+$.
The corollary follows from Proposition \ref{cond-norm}.
\end{proof}

\begin{corollary}\label{coprimos} Suppose $\dim B$ is relatively prime to the rank of $R$.
Then $q: H \to B$ is normal.  \end{corollary}

\begin{proof} We have $\rk R = \dim H_+$. Thus the assumption
implies that $q\vert_{H_+} = \epsilon$ and hence $H_+ \subseteq
H^{\co B}$. Therefore $q$ is normal in view of Corollary
\ref{inclusion}. \end{proof}

The following is an application of Corollary \ref{inclusion} to
the Drinfeld double. See Example \ref{doble}.

\begin{corollary} Let $A$ be a finite dimensional Hopf
algebra and let $q: D(A) \to B$ be a surjective Hopf algebra map.
Suppose that $q\vert_{A} = \epsilon$. Then $q$ is normal. \qed
\end{corollary}

\medbreak Note that $(B, R_q)$ is quasitriangular with $R$-matrix
$R_q = (q\otimes q)(R)$. Note in addition that $f_{R_q} = qf_Rq^*:
B^* \to B$.

Identify $B^*$ with a Hopf subalgebra of $H^*$ by means of the
transpose map $q^*: B^* \hookrightarrow H^*$. Recall  that $H^{\co
B}$ coincides with the space of $B^*$-invariants of $H$ under the
left regular action $\rightharpoonup: B^* \otimes H \to H$, while
${}^{\co B}H$ coincides with the $B^*$-invariants under the right
regular action $\leftharpoonup: H \otimes B^* \to H$.

\begin{lemma}\label{adfR} Let $p \in B^*$, $a \in H^{\co B}$.
Then $a \leftharpoonup p = \ad(f_{R_{21}} p) (a)$. \end{lemma}

\begin{proof} By (QT5) we have $a_2 \otimes a_1 = R^{(1)} a_1 \mathcal S(r^{(1)}) \otimes R^{(2)} a_2
r^{(2)}$, since $R^{-1} = (\mathcal S \otimes \id) (R)$.
Therefore,
\begin{align*} a \leftharpoonup p & = \langle p,  a_1 \rangle a_2 = \langle p,
R^{(2)} a_2 r^{(2)}  \rangle R^{(1)} a_1 \mathcal S(r^{(1)}) \\
& = \langle p_1,  R^{(2)}\rangle \langle p_2, a_2 \rangle \langle
p_3, r^{(2)} \rangle
R^{(1)} a_1 \mathcal S(r^{(1)}) \\
& = f_{R_{21}}(p_1) (p_2 \rightharpoonup a) \mathcal S(f_{R_{21}}p_3) \\
& = f_{R_{21}}(p_1) a \mathcal S(f_{R_{21}}p_2) = f_{R_{21}}(p)_1
a \mathcal S((f_{R_{21}}p)_2) \\ & = \ad(f_{R_{21}} p) (a),
\end{align*}
where we used that $f_{R_{21}}$ is a coalgebra map. This proves
the lemma. \end{proof}

\begin{proposition}\label{centraliz} The following are equivalent:
\begin{itemize}\item[(i)] The map $q:
H \to B$ is normal.  \item[(ii)] $f_{R_{21}}(B^*) \subseteq
(H^{\co B})'$.\end{itemize}
\end{proposition}

Here $(H^{\co B})'$ denotes the centralizer of $H^{\co B}$ in $H$.

\begin{proof} (i) $\Rightarrow$ (ii). By Remark \ref{lcs-rcs} and
Lemma \ref{adfR}, the assumption implies that $\ad(f_{R_{21}}p)
(a) = \epsilon (p) a$, for all $p \in B^*$, $a \in H^{\co B}$.
Then $f_{R_{21}}(B^*) \subseteq (H^{\co B})'$.

(ii) $\Rightarrow$ (i). Let $p \in B^*$, $a \in H^{\co B}$. By
Lemma \ref{adfR}, we have
$$a \leftharpoonup p = \ad(f_{R_{21}}p) (a) = \epsilon (p) a.$$
This shows  that $H^{\co B}$ is invariant under the right regular
action of $B^*$. In view of Remark \ref{lcs-rcs} this implies that
${}^{\co B}H = H^{\co B}$, and the normality of $q$.
\end{proof}

\begin{remark}\label{r-cent} Note that condition (ii) in Proposition \ref{centraliz}
holds in any of the following cases: \begin{itemize}\item[(i)]
$f_R\vert_{B^*} = \epsilon$, or \item[(ii)] $H^{\co B} \subseteq
Z(H)$.\end{itemize}
\end{remark}

\begin{lemma}\label{condic-rk} Let $A \subseteq H^*$ be a Hopf subalgebra.
Suppose that $A \subseteq ({H^*}^{\cop})^{\co f_R}$. Then $A$ is
normal in $H^*$.  \end{lemma}

In particular, if $(\dim A, \rk R) = 1$, then $A$ is normal in
$H^*$. Thus we recover the statement in Corollary \ref{coprimos}.

\begin{proof} Consider the quotient $H \to B = A^*$. The assumption $A \subseteq ({H^*}^{\cop})^{\co
f_R}$ implies that $f_R\vert_{B^*} = \epsilon$. The lemma follows
from Proposition \ref{centraliz}, in view of Remark \ref{r-cent}
(i). \end{proof}

\begin{corollary} Assume $H$ is simple. Then the
restriction  $f_R\vert_{G(H^*)}$ is injective. In particular,
$|G(H^*)|$ divides $\rk R$. \end{corollary}

\begin{proof} Note that the intersection $kG(H^*) \cap ({H^*}^{\cop})^{\co
f_R}$ coincides with the group algebra of the kernel of the
restriction $f_R\vert_{G(H^*)}$. By Lemma \ref{condic-rk}, if $H$
is simple, the restriction $f_R\vert_{G(H^*)}$ must be injective,
as claimed. Moreover, in this case $f_R\vert_{G(H^*)} \simeq
kG(H^*)$ is a Hopf subalgebra of $H_+$ implying that $|G(H^*)|$
divides $\rk R$.
\end{proof}

\begin{proposition} The following are equivalent:
\begin{itemize}\item[(i)] $R_q = 1 \otimes 1$,
\item[(ii)] $f_R (B^*) \subseteq H^{\co q} \cap {}^{\co q}H$,
\item[(iii)] $f_{R_{21}} (B^*) \subseteq H^{\co q} \cap {}^{\co
q}H$. \end{itemize} \end{proposition}

\begin{proof} Clear. See Remark \ref{sub-hopf}. \end{proof}

\section{Hopf algebra quotients arising from $R$-matrices}\label{estructura}

Let $H$ be a finite dimensional Hopf algebra. Then $H$ is
faithfully flat over its left coideal subalgebras \cite{skryabin}.
By \cite[Theorem 3.2]{takeuchi}, the maps
$$I \mapsto H^{\co H/I}, \qquad L \mapsto HL^+,$$ give rise to
inverse bijective correspondences between:
\begin{itemize}\item[(a)] the set of Hopf ideals $I$ of $H$, and
\item[(b)] the set of normal left coideal subalgebras $L$ of $H$.
\end{itemize}

\medbreak \emph{From now on, let $(H, R)$ be a finite dimensional
quasitriangular Hopf algebra. In this section we shall combine the
above with the transmutation theory of Majid \cite[Chapter
9.4]{majid}.}

\medbreak By \cite[Example 9.4.9]{majid} there is a braided Hopf
algebra $\underline H$ in the braided category ${}_H\mathcal M$ of
left $H$-modules, where $\underline H  = H$ as algebras,  with
comultiplication
\begin{equation}\label{delta-braided}\underline \Delta(a) = a_1 \mathcal S(R^{(2)}) \otimes
\ad_{R^{(1)}} (a_2),\end{equation} and $H$ acts on $\underline H$
via the left adjoint action:
$$\ad_h (a) = h_1 a \mathcal S(h_2).$$
Dually \cite[Example 9.4.10]{majid}, there is a braided Hopf
algebra structure $\underline{H^*}$ in $\mathcal M^{H^*} =
{}_H\mathcal M$ where $\underline{H^*} = H^*$ as coalgebras, with
multiplication \begin{equation}\label{m-braided}p \underline{.} q
= \left(\mathcal S(p_1)p_3 \otimes \mathcal S(q_1)\right) (R) \,
p_2q_2,\end{equation} and $H^*$ coacts on $\underline{H^*}$ via
the right adjoint coaction: $\rho: \underline{H^*} \to
\underline{H^*} \otimes H^*$,
\begin{equation}\label{coadjunta}\rho (p) = p_2 \otimes \mathcal
S(p_1)p_3.\end{equation}

In this context, the map $\Phi_R$ given by \eqref{phir} becomes  a
morphism of braided Hopf algebras $\Phi_R: \underline{H^*} \to
\underline{H}$. In other words, $\Phi_R$ is a left $H$-module map
that transforms the product \eqref{m-braided} into the product of
$H$ and the coproduct of $H^*$ into the coproduct
\eqref{delta-braided}. See \cite[Propositions 2.1.14 and
7.4.3]{majid}.

\begin{lemma}\label{left-coid} Let $C \subseteq H^*$ be a subcoalgebra.
Then $\Phi_R(C) \subseteq H$ is a normal left coideal of $H$.
\end{lemma}

\begin{proof} Note that a subcoalgebra $C \subseteq H^*$ is an
$H$-subcomodule under the coadjoint coaction \eqref{coadjunta}.

Because $\Phi_R: \underline{H^*} \to \underline{H}$ is a morphism
of braided Hopf algebras, $\Phi_R$ is a left $H$-module algebra
and coalgebra map.  Therefore, $\Phi_R(C)$ is stable under the
left adjoint action of $H$ (that is, $\Phi_R(C)$ is a normal
subspace of $H$).

Let now $a \in H$ and write $\underline{\Delta} (a) =
\underline{a_1} \otimes \underline{a_2}$. Then we have
\begin{equation}\label{rel-deltas}\Delta(a) = \underline{a_1} R^{(2)} \otimes
\ad_{R^{(1)}}(\underline{a_2}).\end{equation} Indeed, writing $R =
r^{(1)} \otimes r^{(2)} = R^{(1)} \otimes R^{(2)}$, relations
\eqref{inv-r} imply that
$$r^{(1)}R^{(1)} \otimes \mathcal S(R^{(2)}) r^{(2)}
= (\id \otimes \mathcal S) (r^{(1)}R^{(1)} \otimes \mathcal
S^{-1}(r^{(2)})R^{(2)}) = 1 \otimes 1.$$ Identity
\eqref{rel-deltas} follows after combining these with formula
\eqref{delta-braided}.

Since $C \subseteq H^*$ is a subcoalgebra, and $H^* =
\underline{H^*}$ as coalgebras, then $\Phi_R(C) \subseteq
\underline H$ is a subcoalgebra. Identity \eqref{rel-deltas}
implies that $\Delta (\Phi_R(C)) \subseteq H \otimes \Phi_R(C)$
and hence that $\Phi_R(C)$ is a normal left coideal of $H$, as
claimed. \end{proof}

Let $K_C : = k[\Phi_R(C)]$ be the subalgebra of $H$ generated by
$\Phi_R(C)$. By Lemma \ref{left-coid}, $K_C \subseteq H$ is a
normal left coideal subalgebra of $H$.

\begin{theorem}\label{ccociente} Let $C \subseteq H^*$ be a subcoalgebra.
Consider the canonical projection $\pi_C: H \to \overline H_C =
H/HK_C^+$. Then the following hold:
\begin{itemize}\item[(i)] $(\overline H_C, \overline R)$ is a quasitriangular
quotient Hopf algebra with  $\overline R = (\pi_C \otimes
\pi_C)(R)$.

\item[(ii)]$K_C = H^{\co \pi_C}$ and $k[{}_R\Phi(\mathcal S^{-1}C)] = {}^{\co
\pi_C}H$.

\item[(iii)] $\pi_C: H \to \overline H_C$ is the maximal quotient Hopf algebra
with the property \begin{equation}\label{c-radical}(p \otimes
\pi_C)(Q) = p(1)1 , \quad \forall p \in C.\end{equation}
\end{itemize} \end{theorem}

\begin{proof} By \cite{takeuchi} $HK_C^+$ is a Hopf ideal and
$\overline H_C: = H/HK_C^+$ is a quotient Hopf algebra. In
particular, the canonical projection gives a surjective Hopf
algebra map $\pi_C: H \to \overline H_C = H/HK_C^+$, $h \mapsto
\overline h$.

\medbreak Clearly, the quotient $\overline H_C$ is quasitriangular
with $R$-matrix $\overline R = (\pi_C \otimes \pi_C)(R)$. This
proves (i).

\medbreak The results in \cite{takeuchi} imply in addition that
$K_C = H^{\co \pi_C}$. On the other hand, ${}^{\co \pi_C}H =
\mathcal S(H^{\co \pi_C}) = \mathcal S(K_C)$, cf. Section
\ref{apendice}. In view of Lemma \ref{left-right} and the
definition of $K_C$, this implies that ${}^{\co \pi_C}H =
k[{}_R\Phi(\mathcal S^{-1}C)]$. This proves (ii).

\medbreak Next we prove (iii). Let $f \in \overline H_C^*$. Note
that $f(a) = \epsilon(a) f(1)$, for all $a \in H$ such that $a$
belongs to $H^{\co \pi_C}$ or to ${}^{\co \pi_C}H$.

Since $H^{\co \pi_C} = k[\Phi_R(C)]$ is generated by elements
$p(Q^{(1)}) Q^{(2)}$, $p \in C$, then
$$(p \otimes f)(Q) = p(Q^{(1)}) \epsilon(Q^{(2)}) f(1) =
p(1)f(1).$$

Therefore $\pi_C$ satisfies \eqref{c-radical}. To prove
maximality, suppose that $t: H \to B$ is a quotient Hopf algebra
such that $(p \otimes t) (Q) = 1$, for all $p \in C$. Then $t
{\Phi_R}\vert_C = \epsilon 1$, implying that the left coideal
subalgebra $K_C = k[\Phi_R(C)]$ is contained in $H^{\co t}$.
Therefore $\ker \pi_C = HK_C^+ \subseteq H(H^{\co t})^+ = \ker t$
\cite{takeuchi}. Thus $t$ factorizes through $\pi_C$ and
maximality is established. This finishes the proof of the theorem.
\end{proof}

\begin{remark} Suppose $C$ is a one dimensional subcoalgebra; that is, $C = kg$, for some $g \in
G(H^*)$. Lemma \ref{left-coid} implies that $\Phi_R(C)$ is a
normal left coideal of $H$. Since $\dim \Phi_R(C) = 1$,
necessarily $\Phi_R(C) = ka$, where $a \in H$ is a central
group-like element. Thus in this case we recover the result
$\Phi_R(G(H^*)) \subseteq G(H) \cap Z(H)$  in \cite[Theorem 2.3
(b)]{schneider}.
\end{remark}

Theorem \ref{ccociente} applies in particular when $C \subseteq
H^*$ is a Hopf subalgebra. The case $C = H^*$ is contained in the
following theorem. It is of special interest, as it is related
with the notion of \emph{modularization} studied in
\cite{bruguieres, muger}; c.f. Subsection \ref{modular}.

\begin{theorem}\label{cociente} Let $K = \Phi_R(H^*) \subseteq H$. Then $K$ is a
normal left coideal subalgebra of $H$.

Consider the canonical projection $\pi: H \to \overline H =
H/HK^+$. Then the following hold:
\begin{itemize}\item[(i)] $(\overline H, \overline R)$ is a triangular
quotient Hopf algebra with $R$-matrix $\overline R = (\pi \otimes
\pi)(R)$.

\item[(ii)]$\Phi_R(H^*) = H^{\co \pi}$ and ${}_R\Phi(H^*) = {}^{\co
\pi}H$.

\item[(iii)] $\pi: H \to \overline H$ is the maximal quotient Hopf
algebra with the property \begin{equation}\label{radical}(\pi
\otimes \id)(Q) = 1 = (\id \otimes \pi)(Q).\end{equation}
\end{itemize} \end{theorem}

Note that finite dimensional triangular Hopf algebras  have been
completely classified \cite{eg2, aeg, eg-triangular}. It turns out
that $(\overline H, \overline R)$ is necessarily a twisting of a
modified supergroup algebra. We shall use this fact in later
applications of this theorem. See Section \ref{aplicaciones}.

\begin{proof} Parts (ii) and (iii) are special instances of Theorem
\ref{ccociente}. Indeed, $\Phi_R(H^*)$ is a subalgebra of $H$. To
prove part (i) it is enough to show that $\overline H$ is
triangular. We have $\Phi_{\overline R} = \pi \Phi_R \pi^*$, where
$\pi^*: \overline{H}^* \hookrightarrow H^*$ is the dual inclusion.
Hence $\Phi_{\overline R} = \epsilon 1$, because $\pi\vert_K =
\epsilon 1$. Therefore $(\overline H, \overline R)$ is triangular,
as claimed.
\end{proof}

\begin{remark}\label{rs-min} Let $C \subseteq H^*$ be a
subcoalgebra. Keep the notation in Theorems \ref{ccociente},
\ref{cociente}.

(1) Since $\Phi_R(H^*) \subseteq H_R$ \cite{radford}, then the
left coideal $\Phi_R(C)$ and also $K_C = k[\Phi_R(C)]$ are
contained in $H_R$.
 Moreover, since $Q \in H_R \otimes H_R$, then
$\Phi_R(H^*) = \Phi_R(H_R^*)$. In particular,  whenever $H_R$ is
factorizable, $\Phi_R(H^*) = H_R$ is a normal Hopf subalgebra of
$H$.

(2) The inclusion $K_C \subseteq H_R$ also implies that $[H: H_R]$
divides $\dim \overline H_C$.

Indeed, consider the exact sequences $K_C \subseteq H \to
\overline H_C$, $K_C \subseteq H_R \to (\overline{H_R})_C$. Taking
dimensions we get
\begin{align*}\dim H & = \dim H_R [H : H_R] = \dim K_C \dim
(\overline{H_R})_C \, [H: H_R] \\ & = \dim K_C \dim \overline
H_C.\end{align*} Hence $\dim \overline H_C = \dim
(\overline{H_R})_C \, [H: H_R]$, and $[H : H_R] / \dim \overline
H_C$, as claimed.

(3) The maximality of the quotient $H \to \overline H_C$ implies
that there is a sequence of surjective Hopf algebra maps $H \to
\overline H_C \to \overline H$.
\end{remark}

\begin{lemma} Let $C \subseteq H^*$ be a subcoalgebra.

\begin{itemize}\item[(i)] $\overline H_C = H$ if and only if $\Phi_R\vert_C =
\epsilon$ if and only if $C \subseteq \overline H^*$.

\item[(ii)] If $\overline H_C = k$ then $H$ is factorizable.

\item[(iii)] $k[C] \supseteq \overline H_C^*$ if and only if $\overline
H_C$ is triangular. \end{itemize}
\end{lemma}

In particular, if $H$ is not factorizable and $C \nsubseteq
\overline H^*$, then $\overline H_C$ is a proper quotient Hopf
algebra.

\begin{proof} Part (i) follows from the definition of $\overline
H_C$. Part (ii) follows from Remark \ref{rs-min} (3). Finally,
part (iii) is a consequence of Theorem \ref{ccociente} (iii).
\end{proof}

\begin{theorem}\label{indice} Suppose $A \subseteq H^*$ is a Hopf subalgebra. Then there is a
Hopf subalgebra $B \subseteq H^*$ such that

\begin{itemize}\item[(i)] $[H^*: A]$ divides $\dim B$.

\item[(ii)] $\Phi_R\vert_{A\cap B} = \epsilon$. Furthermore, $B =
H^*$ if and only if $\Phi_R\vert_{A} = \epsilon$. \end{itemize}
\end{theorem}

\begin{proof} Apply Theorem \ref{ccociente} to the subcoalgebra
$A$. Then $\Phi_R(A)$ is a subalgebra of $H$ and thus $K_A =
\Phi_R(A)$.

Moreover, $\Phi_R: A \to \Phi_R(A)$ is a surjective map of braided
Hopf algebras over $H$. Hence $\dim \Phi_R(A) / \dim A$ by the
Nichols--Zoeller theorem applied to the corresponding biproducts.

Consider the quotient $\overline H_C = H/HK_C^+$. So that $\dim
\overline H_C = \frac{\dim H}{\dim \Phi_R(A)}$. Let $B = \overline
H_A^* \subseteq H^*$. Then $B$ satisfies (i) and (ii).
\end{proof}

The following application of Theorem \ref{ccociente} gives
restrictions on the possible quotients of a factorizable Hopf
algebra.

\begin{theorem} Let $(H, R)$ be a factorizable Hopf algebra.
Suppose $A \subseteq H^*$ is a Hopf subalgebra. Then there is a
Hopf subalgebra $B \subseteq H^*$ such that

\begin{itemize}\item[(i)] $\dim A \dim B = \dim H$.

\item[(ii)] $A\cap B = k1$.\end{itemize} \end{theorem}

\begin{proof} In this case $\dim \Phi_R(A) = \dim A$ because $H$ is factorizable.
As in the proof of Theorem \ref{indice}, take for $B$ the Hopf
subalgebra $\overline H_A^* \subseteq H^*$. Then $B$ satisfies
(i). Moreover, we have $\Phi_R\vert_{A\cap B} = \epsilon$, whence
$A\cap B = k1$, since $H$ is factorizable by assumption. Thus $B$
satisfies also (ii).
\end{proof}

Part (ii) of the following proposition generalizes \cite[Theorem
2.3 (b)]{schneider}. Part (iii) gives a refinement of
\cite[Proposition 3 c)]{radford}.

\begin{proposition}\label{central-gl} Identify $\overline H^*$ with a Hopf subalgebra of $H^*$. Then
\begin{itemize}\item[(i)]$\Phi_R$ induces an isomorphism $\Phi_R: H^*/(\overline H^*)^+H^* \overset{\simeq}\to K$.
\item[(ii)] $\Phi_R$ induces an injective group homomorphism $$\Phi_R: G(H^*)/G(\overline H^*) \hookrightarrow G(H) \cap Z(H).$$
\item[(iii)] Let $\alpha \in G(H^*)$ be the modular element. Then $\alpha \in
G(\overline H^*)$. \end{itemize}

In particular, if $G(\overline H^*) = 1$, then $H$ is unimodular.
\end{proposition}

\begin{proof} (i). We shall see that $\ker \Phi_R = (\overline H^*)^+H^*$, whence the claimed isomorphism.
We know that $$\dim K = \dfrac{\dim H}{\dim \overline H} = \dim
(H^*/(\overline H^*)^+H^*) = \dim H^* - \dim (\overline
H^*)^+H^*.$$ Thus it will be enough to see that $(\overline
H^*)^+H^* \subseteq \ker \Phi_R$. For this, let $p \in (\overline
H^*)^+$, $f \in H^*$. So that $\Phi_R(p) = \epsilon(p) = 0$. Using
(QT1), (QT3), we have
\begin{align*} \phi_R(pf) & = (pf)(r^{(2)}R^{(1)}) r^{(1)}R^{(2)} \\
& = (p \otimes f)(\Delta(r^{(2)})\Delta(R^{(1)})) r^{(1)}R^{(2)} \\
& = (p \otimes f) (s^{(2)}R^{(1)} \otimes r^{(2)}t^{(1)}) r^{(1)}s^{(1)}R^{(2)}t^{(2)} \\
& = f(r^{(2)}t^{(1)}) r^{(1)}\Phi_R(p) t^{(2)} = 0,
\end{align*}
where $R = r = s = t$. This gives the desired inclusion.

(ii). The irreducible character $\chi$ of $H$ comes from an
irreducible character of $\overline H$ if and only if
$\Phi_R(\chi) = \deg \chi$. Thus $\ker \Phi_R\vert_{G(H^*)} =
G(\overline H^*)$.

(iii). By \cite{radford} we have $\Phi_R(\alpha) = 1$. Using (ii),
we get $\alpha \in G(\overline H^*)$.
\end{proof}

\medbreak We next give some sufficient conditions for the map
$\pi$ to be normal.

\begin{proposition}\label{normalidad} Suppose that either one of the following conditions hold:

\begin{itemize}\item[(a)] $H_+ \subseteq K$,

\item[(b)] $R_{21}R = RR_{21}$,

\item[(c)] $\overline R = 1 \otimes 1$ and $K$ is commutative, or

\item[(d)] $\overline R = 1 \otimes 1$ and $\dim K$ is relatively prime
to $\dim \overline H$. \end{itemize}

Then $K \subseteq H$ is a normal Hopf subalgebra and there is an
exact sequence of Hopf algebras
$$k \to K \hookrightarrow H  \overset{\pi}\to \overline H \to k.$$ \end{proposition}

\begin{proof} Suppose (a) holds. Then the conclusion follows from Corollary \ref{inclusion}.

Now assume (b). Then $H^{\co \pi} = {}^{\co \pi}H$, in view of
Theorem \ref{cociente} (ii). Hence $\pi$ is normal as claimed.

If (c) holds, then $f_R(\overline H^*) \subseteq K = H^{\co \pi}$,
and the conclusion follows from Proposition \ref{centraliz}.

Finally, in case (d), necessarily $f_R\vert_{\overline H^*} =
\epsilon$, and the conclusion follows from Remark \ref{r-cent}.
\end{proof}

\begin{theorem}\label{hr-coconm} Suppose that $R \in kG(H) \otimes kG(H)$. Then
$\Phi_R(H^*) \subseteq H$ is a commutative normal Hopf subalgebra,
and there is an exact sequence of Hopf algebras $$k \to
\Phi_R(H^*) \to H \overset{\pi}\to \overline H \to k.$$ In
particular, $H$ is an extension of a dual group algebra by a
twisting of a modified supergroup algebra.
\end{theorem}

\begin{proof} It will be enough to show that $\Phi_R(H^*)$ is a Hopf
subalgebra. Commutativity of $\Phi_R(H^*) \subseteq (kG(H))_R$
will follow from Example \ref{cocomm}. The assumption $R \in kG(H)
\otimes kG(H)$ is equivalent to $H_R \subseteq kG(H)$. If this
holds, then $\Phi_R(H^*) \subseteq H_R \subseteq kG(H)$, and
therefore $\Phi_R(H^*)$ is necessarily a Hopf subalgebra, as
claimed. \end{proof}

\begin{example} Quasitriangular structures on nonsemisimple Hopf
algebras of dimension $8$ have been classified in \cite{wakui}. It
turns out that there is only one example, denoted $A_{C_2}$, of
such a Hopf algebra which admits an $R$-matrix that is not
triangular.

Explicitly, $H = A_{C_2}$ is presented by generators $g, x, y$ and
relations $g^2 = 1$, $x^2 = y^2 = 0$, $gx = -xg$, $gy = -yg$, $xy
= -yx$, with $\Delta(g) = g \otimes g$, $\Delta(x) = x \otimes g +
1 \otimes x$, $\Delta(y) = y \otimes g + 1 \otimes y$. The
$R$-matrices on $H$ are parameterized by $R_{abcd}$, with $a, b,
c, d \in k$, such that $R_{abcd}$ is minimal whenever $b \neq 0$
or $c \neq 0$ or $ad \neq 0$.

By \cite[Proof of Proposition 2.9]{wakui}, for $R = R_{abcd}$, we
have $$R_{21}R = 1 \otimes 1 + (b-c) (y \otimes x - x \otimes y)
(1 \otimes g) - (b -c)^2 xy \otimes xy.$$ Therefore, $H$ is
triangular if $b = c$, and otherwise $\Phi_R(H^*)$ is a left
coideal subalgebra of dimension  $4$  spanned by $1, xg, yg, xy$.
In particular, $\Phi_R(H^*)$ is not a Hopf subalgebra of $H$.
\end{example}

\begin{example}\label{gelaki's} Let $p \neq q$ be prime numbers such that $p = 1$ (mod $q$).
Let $A_i$, $0 \leq i \leq p-1$, be one of the nontrivial self-dual
semisimple Hopf algebras  with $\dim A_i = pq^2$ studied in
\cite{gelaki, examples}. Quasitriangular structures in $A_i$'s
were discussed in \cite{gelaki}. We list some properties of these
Hopf algebras (see {\it loc. cit.} for a proof):

\begin{itemize} \item[(i)] $G(A_0) \simeq \mathbb Z_q \times \mathbb Z_q$ and $G(A_i) \simeq \mathbb Z_{q^2}$, $i >0$.

\item[(ii)] $G(A_i) \cap Z(A_i) \simeq \mathbb Z_q$, and this is the only
normal Hopf algebra quotient of dimension $q$ of $A_i^*$.

\item[(iii)] Suppose $i > 0$. Then $A_i$ is quasitriangular if and only if  $q = 2$.
Moreover, when $q = 2$, $A$ admits $2p-2$ minimal quasitriangular structures, none of them triangular.

\item[(iv)] $A_0$ admits triangular structure with $(A_0)_R = kG(A_0)$. $A_0$ is minimal quasitriangular if
and only if $q = 2$ and in this case none of the quasitriangular structures is triangular.

\item[(v)] The proper Hopf subalgebras of $A_i$ are commutative or cocommutative. In particular,
if $A_i$ admits a quasitriangular structure which is not minimal, then $(A_i)_R \subseteq kG(A_i)$.

Furthermore, $A_i$ has a unique normal quotient of dimension $q$;
thus a unique Hopf subalgebra of  dimension $pq$ (which is
necessarily normal, because $q$ is the smallest prime number
dividing $\dim A$). Hence if $(A_i, R)$ is minimal
quasitriangular, necessarily $(A_i)_+ = (A_i)_- = A$. That is,
$f_R: A_i^{*\cop} \to A_i$ is an isomorphism.
\end{itemize}

We prove next the following additional properties.

\begin{lemma}\label{ej-gel}  Let $A = A_i$, $0 \leq i \leq p-1$, and assume $(A, R)$ is quasitriangular. Then we have
\begin{itemize}\item[(vi)] $A$ is not factorizable.
\end{itemize}
\noindent Suppose $q$ is odd and $(A, R)$ is not triangular. Then
\begin{itemize}
\item[(vii)] $\Phi_R(A^*)$ is a commutative (normal) Hopf subalgebra of dimension $q$ or $pq$.

\item[(viii)] $G(A) \cap Z(A) \subseteq \Phi_R(A^*)$.
\end{itemize} \end{lemma}

\begin{proof} (vi). Suppose on the contrary that $A$ is factorizable.
By Proposition \ref{central-gl}, $|G(\overline{A}^*)| = q$. This
implies that $\overline A \neq 1$ and hence that $A$ is not
factorizable, thus proving (vi).

\medbreak Now suppose $q$ is odd and $(A, R)$ is not triangular.
Let $K = \Phi_R(A^*)$  and consider the sequence $K \subseteq A
\to \overline A$. By (vi),  $\overline A \neq 1$, therefore $\dim
\overline A = q^2, pq$ or $q$. By (iii) and (iv) $A$ is not
minimal. Then by (v) and Theorem \ref{hr-coconm} $K$ is a
(commutative) normal Hopf subalgebra of $A$. Necessarily $\dim K =
q$ or $pq$, thus part (vii) follows.

We have $\overline H \simeq kF$ for some group $F$ with $|F| = q$
or $pq$. In particular, $G(k^F)$, being a subgroup of $G(A^*)
\simeq \mathbb Z_q \times \mathbb Z_q$, is of order $q$.

By Proposition \ref{central-gl} (i), $\Phi_R$ induces an injective group homomorphism $G(A^*)/G(k^F) \to G(A) \cap Z(A)$. Since $G(A^*)/G(k^F)$ is not trivial and $G(A) \cap Z(A)$ is of order $q$, this homomorphism is in fact an isomorphism $G(A^*)/G(k^F) \simeq G(A) \cap Z(A)$. This proves (viii).
\end{proof}

\end{example}

\subsection{Minimality} The following lemma gives a relation between the minimal
quasitriangular Hopf subalgebras of $H$ and $\overline H$.

\begin{lemma}If $\overline H$ is minimal, then $H$ is minimal. \end{lemma}

\begin{proof} Consider the sequence $H_R \subseteq H \overset{\pi}\to \overline
H$. Since $\Phi_R(H^*) \subseteq H_R$,
$$H_R^{\co \pi} = H_R \cap H^{\co \pi} = H_R \cap \Phi_R(H^*) = \Phi_R(H^*) = \Phi_R(H_R^*).$$
Therefore $\pi(H_R) \simeq H_R/\Phi_R(H_R^*)^+H_R =
\overline{H_R}$. Hence we have exact sequences $$\Phi_R(H^*)
\subseteq H \to \overline H, \quad \Phi_R(H^*) \subseteq H_R \to
\overline{H_R}.$$ Taking dimensions we get $$\dim H = \dim
\Phi_R(H^*) \dim \overline H, \quad \dim H_R = \dim \Phi_R(H^*)
\dim \pi(H_R).$$ Since $\overline R \in \pi(H_R) \otimes
\pi(H_R)$, then $\overline H_{\overline R} \subseteq \pi(H_R)$.
Now assume that $\overline H$ is minimal; then $\overline H =
\overline H_{\overline R} = \pi(H_R)$. Therefore $\dim H = \dim
H_R$ and $H = H_R$ is minimal, as claimed.
\end{proof}

\subsection{Modularization (\cite{bruguieres, muger}.)}\label{modular} The quotient $\pi: H \to \overline H$ is an
analogue of the \emph{transparent} tensor subcategory $\mathcal
T_{\mathcal C}$ of the category  $\mathcal C = \Rep H$ \cite[pp.
224]{bruguieres}.

Indeed, if $H$ is semisimple, then this follows from relations in
Theorem \ref{cociente} (iii). However, our construction does not
assume semisimplicity nor even existence of a ribbon structure.

\smallbreak \emph{Assume  that $H$ is semisimple.} Then the
category $H$-mod is a {\it premodular category}; here the braiding
is given by the action of the $R$-matrix:
$$\sigma_{V, W} (v \otimes w) = R^{(2)}w \otimes R^{(1)}v, \quad v \in V, \, w \in W.$$
The ribbon structure is given by the action of the (central)
Drinfeld element $u = \mathcal S (R^{(2)}) R^{(1)}$.

\medbreak Recall from \cite{bruguieres} that a left $H$-module $V$
is called {\it transparent} if for any  $H$-module $W$,
$\sigma_{W, V} \sigma_{V, W} = \id$. Observe that $(H, R)$ is
triangular if and only if all simple $H$-modules are transparent,
and $(H, R)$ is factorizable if and only if $H$ has no non-trivial
irreducible transparent simple modules.

Let $R(H) \subseteq H^*$ be the character algebra of $H$.
Condition (QT5) implies that $R(H)$ is a commutative semisimple algebra over
$k$.
Let $\mathcal T (H) \subset R(H)$ denote the linear span of the characters of
the irreducible transparent $H$-modules.

\begin{proposition} $\mathcal T (H) = R(\overline H) \subseteq
R(H)$. \end{proposition}

\begin{proof} The irreducible character $\chi \in
R(H)$ is transparent if and only if $\Phi_R(\chi) = \chi(1) 1$.
The proof follows from Theorem \ref{cociente} (iii).
\end{proof}

By \cite{eg}, after modifying if necessary the $R$-matrix
$\overline R$, $\overline H$ is isomorphic to a twisting of a
group algebra $H \simeq (kG)_J$, where $G$ is a finite group and
$J \in kG \otimes kG$ is a normalized $2$-cocycle.

\begin{remark} Note that $\Rep H$, with its canonical braiding, is modularizable in the sense of
\cite{bruguieres} if and only $\overline u = 1$, where $\overline
u$ is the Drinfeld element of $\overline H$, that coincides with
the image of the Drinfeld element $u \in H$  under the canonical
projection $H \to \overline H$. Indeed the last condition amounts
to the category of transparent objects in $\Rep H$ being tannakian
\cite{deligne, eg}.
\end{remark}

\section{Application}\label{aplicaciones}

In this section we apply the results in Sections \ref{apendice}
and \ref{estructura} to the classification of finite dimensional
Hopf algebras with square-free dimension.

\smallbreak In what follows we shall consider  a quasitriangular
Hopf algebra $(H, R)$. Assume that $\dim H$ is odd and
square-free.

\begin{lemma}\label{triang-ss} Suppose $(H, R)$ is triangular. Then $H$ is
semisimple and $R = 1 \otimes 1$. In particular, $H$ is isomorphic
to a group algebra. \end{lemma}

\begin{proof} It follows from \cite[Theorem 4.3]{eg-triangular} that $H$ has the Chevalley property;
in other words, the coradical of the dual Hopf algebra $H^*$ is a
Hopf subalgebra. By \cite{aeg}, $H$ is twist equivalent to a
triangular Hopf algebra $H'$ with $R$-matrix $R'$ of rank $\leq
2$.

Let $u' \in {H'}_{R'}$ be the Drinfeld element. Since $H'$ is
triangular, $u' \in G(H')$. Moreover, ${u'}^2 = 1$, since $\rk R'
\leq 2$. Because $\dim H' = \dim H$ is odd, we get $u' = 1$. Hence
$\mathcal S^2 = \ad_{u'} = \id$, implying that $H'$, and therefore
also $H$, are semisimple.

Because $H$ is odd dimensional, its Drinfeld element $u$ is
necessarily trivial (otherwise it would be a group-like element of
order $2$, by \cite[Lemma 2.1]{eg}). By \cite[Corollary
2.2.2]{eg-reptriang}, $R = 1 \otimes 1$, since $\dim H_R$ must be
a square dividing $\dim H$.
\end{proof}

\begin{lemma}\label{sf-ext} We have $\Phi_R(H^*) =
H_R$ is a normal Hopf subalgebra of $H$ and there is an exact
sequence of Hopf algebras $$k \to H_R \to H \overset{\pi}\to
\overline H \to k.$$
\end{lemma}

\begin{proof} Consider the projection $\pi: H \to \overline H$ and let $K = H^{\co \pi} = \Phi_R(H^*)$.
By Lemma \ref{triang-ss}, $\overline R = 1 \otimes 1$ and
therefore $f_R(\overline H^*) \subseteq K$.

Because $\dim H$ is square-free, $\dim \overline H$ and $\dim
H^{\co \pi}$ are relatively prime.  This implies that
$f_R(\overline H^*) = k1$ since,  by the Nichols-Zoeller Theorem,
$\dim f_R(\overline H^*)$ divides both $\dim \overline H$ and
$\dim H^{\co \pi}$. In view of Proposition \ref{centraliz}, $\pi$
is normal (c. f. Remark \ref{r-cent}).

\smallbreak It remains to show that $K = H_R$ or, in other words,
that $H_R$ is factorizable. This is in turn equivalent to proving
that $\overline{H_R} = k$.  Since $\dim H_R$ is square-free and
$\dim H_R / (\dim H_+)^2$, we have $H_+ = H_- = H_R$. In
particular, $R$ induces an isomorphism $f_R: H_R^{*\cop} \to H_R$.

Since $\overline R = 1 \otimes 1$,  as in the previous paragraph,
we get $f_R(\overline{H_R}^*) = k1$. This implies that
$\overline{H_R} = k$ as claimed. The proof of the lemma is now
complete.
\end{proof}

\begin{proposition} $H$ is semisimple. \end{proposition}

\begin{proof} If $(H, R)$ is triangular, the proposition follows from
Lemma \ref{triang-ss}. Suppose next that $(H, R)$ is factorizable.
Then $H$ is minimal quasitriangular and therefore, $\dim H$ being
square-free, $\dim H = \rk R$. Thus $H$ is semisimple, by Lemma
\ref{fact-odd}. The general case follows from Lemma \ref{sf-ext},
since an extension of semisimple Hopf algebras is semisimple.
\end{proof}

\subsection*{Proof of Theorem \ref{ap1}} By \cite[Theorem 3.2]{schneider} the square of the dimensions of the
irreducible $H_R$-modules divide $\dim H_R$. Since $\dim H_R$ is
square-free all irreducible $H_R$-modules must be one-dimensional
and therefore $H_R$ is commutative.

Hence $H_R \simeq k^{\Gamma}$ for some finite group $\Gamma$. In
addition $\Gamma$ must be abelian (hence cyclic) because
$k^{\Gamma}$ is quasitriangular.

\medbreak In view of Lemma \ref{sf-ext} we have an \emph{abelian}
extension \cite{Maext}
$$k \to k^{\Gamma} \to H \to kF \to k,$$ where $\Gamma$ is cyclic
and $R \in k^{\Gamma} \otimes k^{\Gamma}$. This extension gives
rise to compatible actions $\triangleright: \Gamma \times F \to F$
and $\triangleleft: \Gamma \times F \to \Gamma$ and 2-cocycles
$\sigma: F \times F \to (k^{\Gamma})^{\times}$, $\tau: \Gamma
\times \Gamma \to (k^F)^{\times}$ in such a way that $H$ is a
bicrossed product $H \simeq k^{\Gamma} {}^{\tau}\#_{\sigma} kF$.
Moreover, since $\Gamma$ is cyclic we may assume (up to Hopf
algebra isomorphisms) that $\tau = 1$.

\medbreak In the basis $(e_sx: \, s \in \Gamma, \, x \in F)$, the
comultiplication and multiplication of $k^{\Gamma}
{}^{\tau}\#_{\sigma} kF$ are determined by
\begin{equation}\label{delta}\Delta(e_gx) = \sum_{st = g} e_s \# (t \triangleright
x) \otimes e_t \# x,\end{equation}
$$(e_s \# x)(e_t \# y) = \delta_{s \triangleleft x, t} e_s\sigma(x, y) \# xy.$$

On the other hand the $R$-matrix $R \in k^{\Gamma} \otimes
k^{\Gamma}$ can be written in the form $R = \sum_{s, t} \langle s,
t \rangle e_s \otimes e_t$, for some (non-degenerate) bicharacter
$\langle \, , \, \rangle : \Gamma \times \Gamma \to k^{\times}$.

Let $x \in F$. We have $\Delta^{\cop}(1 \# x) =  \sum_{s, t}  e_t
\# x \otimes e_s \# (t \triangleright x)$.

On the other hand,
$$R\Delta(1 \# x)R^{-1} = \sum_{s, t} \langle s, t \rangle
\langle s \triangleleft (t \triangleright x), t \triangleleft x
\rangle^{-1} e_s \# (t \triangleright x) \otimes e_t \# x.$$
Comparing both expressions we find that the action
$\triangleright: \Gamma \times F \to F$ must be trivial. This
implies that $H$ is cocommutative in view of \eqref{delta}.

\medbreak So $H \simeq kG$, and $G$ fits into an extension $1 \to
\widehat \Gamma \to G \to F \to 1$. Then $G \simeq \widehat \Gamma
\rtimes F$ by the Schur-Zassenhaus' Lemma.

Also $(g^{-1} \otimes g^{-1})R(g \otimes g) = R$, for all $g \in
G$, implying that the form $\langle \, , \, \rangle$ is
$F$-invariant.
 This finishes the proof of the theorem. \qed

\bibliographystyle{amsalpha}

\end{document}